\documentclass[%
 aip,
 amsmath,amssymb,
 reprint,%
]{revtex4-1}

\usepackage{graphicx}
\usepackage{dcolumn}
\usepackage{bm}

\usepackage[utf8]{inputenc}
\usepackage[T1]{fontenc}
\usepackage{mathptmx}
\usepackage{etoolbox}

\makeatletter
\def\@email#1#2{%
 \endgroup
 \patchcmd{\titleblock@produce}
  {\frontmatter@RRAPformat}
  {\frontmatter@RRAPformat{\produce@RRAP{*#1\href{mailto:#2}{#2}}}\frontmatter@RRAPformat}
  {}{}
}%
\makeatother
\begin{document}


\title[Stiff Neural Ordinary Differential Equations]{Stiff Neural Ordinary Differential Equations}
\author{Suyong Kim}
\affiliation{Department of Mechanical Engineering, Massachusetts Institute of Technology, Cambridge, MA 02139, USA}

\author{Weiqi Ji}
\affiliation{Department of Mechanical Engineering, Massachusetts Institute of Technology, Cambridge, MA 02139, USA}

\author{Sili Deng}
\thanks{Correspond to silideng@mit.edu (SD)}
\affiliation{Department of Mechanical Engineering, Massachusetts Institute of Technology, Cambridge, MA 02139, USA}

\author{Yingbo Ma}
\affiliation{Julia Computing Inc., Cambridge, MA 02144, USA}

\author{Christopher Rackauckas}
\thanks{Correspond to contact@chrisrackauckas.com (CR)}
\affiliation{Department of Mathematics, Massachusetts Institute of Technology, Cambridge, MA 02139, USA}
\affiliation{School of Pharmacy, University of Maryland, Baltimore, MD 21201, USA}
\affiliation{Pumas AI, Baltimore, MD 21201, USA}

\date{\today}

\begin{abstract}
Neural Ordinary Differential Equations (ODEs) are a promising approach to learn dynamical models from time-series data in science and engineering applications. This work aims at learning neural ODE for stiff systems, which are usually raised from chemical kinetic modeling in chemical and biological systems. We first show the challenges of learning neural ODEs in the classical stiff ODE systems of Robertson’s problem and propose techniques to mitigate the challenges associated with scale separations in stiff systems. We then present successful demonstrations in stiff systems of Robertson’s problem and an air pollution problem. The demonstrations show that the usage of deep networks with rectified activations, proper scaling of the network outputs as well as loss functions, and stabilized gradient calculations are the key techniques enabling the learning of stiff neural ODEs. The success of learning stiff neural ODEs opens up possibilities of using neural ODEs in applications with widely varying time-scales, like chemical dynamics in energy conversion, environmental engineering, and life sciences. 
\end{abstract}

\maketitle

\begin{quotation}
Neural Ordinary Differential Equations (ODEs) have been emerged as a powerful tool to describe a dynamical system using an artificial neural network. Despite its many advantages, there are many scientific and engineering examples where neural ODEs may fail due to stiffness. This study demonstrates how instability arises in neural ODEs during the training with the benchmark stiff system. Further, we apply scaling to differential equations and loss functions to mitigate stiffness. The proposed technique in adjoint handling and equation scaling for stiff neural ODEs can be used in many biological, chemical, and environmental systems subject to stiffness.

\end{quotation}

\section{\label{sec:intro}Introduction}

Neural Ordinary Differential Equations (ODEs) are elegant approaches for data-driven modeling of dynamical systems from time-series data in science and engineering applications. Neural ODEs offer various advantages over physics-based modeling and traditional neural network modeling. Compared to physics-based modeling where a presumed functional form is required, the neural ODEs approach precludes the need for expert knowledge in identifying such prior information \cite{rackauckas2020universal,ji2020autonomous}. Compared to traditional neural network modeling, such as recurrent neural networks, neural ODEs enable flexibility in the modeling of irregularly and incomplete sampled time series \cite{rubanova2019latent}, and can be more efficient by taking the advantages of modern ODE solvers and adjoint methods \cite{chen2018neural,rackauckas2019diffeqflux}.

Despite the success of neural ODEs \cite{chen2018neural} in many sciences and engineering problems \cite{bills2020universal,portwood2019turbulence,maulik2020time,owoyele2020neural}, it has been recognized that it is very challenging to learn neural ODEs for stiff dynamics \cite{ghosh2020steer}, characterized by dynamics with widely separated time scales. Part of the challenge is due to the high computational cost of solving stiff ODEs as well as handling the ill-conditioned gradients of loss functions with respect to neural network weights \cite{anantharaman2020accelerating}. In addition, it has been shown that stiffness could lead to failures in many data-driven modeling approaches, such as reduced order of modeling \cite{huang2020model} and physics-informed neural networks \cite{ji2020stiff}. It was suggested that the stiffness could lead to gradient pathologies and ill-conditioned optimization problems, which leads to the failure of stochastic gradient descent based optimization \cite{wang2020understanding}.

This work is motivated to elucidate the challenges of learning neural ODEs for stiff systems and proposes strategies to overcome these obstacles. We study learning neural ODEs on data generated from two classical stiff systems, ROBER \cite{robertson1966solution} and POLLU \cite{verwer1994gauss}, which are extensively used for benchmarking stiff ODE integrators. Both ROBER and POLLU describe the dynamics of species concentrations in stiff chemical reaction systems. We propose strategies to mitigate the challenges associated with the scale separations in time as well as the magnitudes of states. We are thus able to successfully learn stiff neural ODEs for those two benchmark problems. We demonstrate how the application of adjoint methods to stiff systems leads to numerical issues and changes in the computational complexity that are overcome by proposing new adjoint techniques. Those strategies not only demonstrate successful numerical methods and architectures for handling time scale separations in neural ODEs but also elucidate general theoretical issues of data-driven modeling of stiff dynamical systems which can be transferred to other techniques.

\section{Background}
\label{sec:background}

\subsection{Neural Ordinary Differential Equations}

Many science and engineering problems can be formulated as ODEs, that is,

\begin{equation}
\frac{dy\left(t \right)}{dt} = f \left(y \left(t \right), \theta, t \right) \label{eq:eq1_1}
\end{equation}

where $t$ is time, $y(t)$ is the state variables, and the function $f$ models the dynamics. Identifying the model $f$ is one of the central tasks in many scientific discoveries, which usually involves proposing functional forms and then performing parameter inference against observation data. Proposing functional forms of $f$ is a very challenging task for many complex systems, such as modeling the gene regulatory networks and cell signaling networks in life science \cite{hoffmann2019reactive, yuan2020cellbox}. Discovering the functional forms of $f$ could require decades of effort with expert knowledge and it is often the case that a lot of unknown dynamics are yet to discover in those complex systems \cite{brunton2016discovering,mangan2016inferring}.
With the help of the universal approximation theorem, one can approximate the model $f$ using a neural network without worrying about missing important dynamics as in physics-based modeling, i.e.,
\begin{equation}
\frac{dy\left(t \right)}{dt} = NN_{\theta} \left(y \left(t \right) \right) \label{eq:eq1_2}
\end{equation}

If one can measure the data pair of $\left(y\left(t\right), \frac{dy\left(t\right)}{dt}\right)$, one can train the neural network straightforwardly as a normal supervised learning problem. However, we usually only have access to the time-series data of $y\left(t\right)$. Instead, neural ODEs train the $NN$ by integrating the ODEs,
\begin{equation}
\begin{split}
y\left(t_1\right) &=y\left(t_{0}\right)+\int_{t_{0}}^{t_1} NN \left(y \left(t \right) \right) dt \\
&= ODESolve\left(y\left(t_{0}\right),NN,t_{0},t_{1},\theta\right)\label{eq:eq1_3}
\end{split}
\end{equation}
in which $y(t_0)$ are the initial conditions of the integration, $t \in (t_0,t_1)$ is the range of the integration. We can define the loss functions as the difference between the observed state variables and the predicted state variables,
\begin{equation}
L(\theta) = MAE \left(y \left(t \right)^{model}, y \left(t \right)^{obs} \right) \label{eq:eq1_4}
\end{equation}
where the mean absolute error (MAE) as an example metric. Using modern differentiable ODE solvers \cite{rackauckas2019diffeqflux}, one can compute the gradient of the loss function to the model parameters efficiently. This work employs Julia’s~\cite{Julia-2017} Scientific Machine Learning ecosystem, mainly consisting of DifferentialEquations.jl \cite{rackauckas2017differentialequations}, Flux.jl \cite{innes2018flux}, ForwardDiff.jl \cite{revels2016forward}, to integrate the neural ODEs and do the backpropagation. We can thus use stochastic gradient descent algorithms to optimize the neural network.

\subsection{Stiff ODE Systems}

Ernst Hairer's classic working definition for the field of numerical ODE solving is ``stiff equations are problems for which explicit methods don't work''~\cite{wanner1996solving,shampine1979user}. In other words, it is the case where numerical ODE solvers like the \verb dopri  method, chosen in the original neural ODE paper \cite{chen2018neural}, fail to adequately solve the differential equation. Many definitions have been proposed to explain this phenomena. A commonly used definition is the stiffness index:
\begin{equation}
    S = \frac{Re(\lambda_{max})}{Re(\lambda_{min})} (t_{1} - t_0)
\end{equation}
where $\lambda_i$ are the eigenvalues of the Jacobian. While these types of stiffness can be helpful in identifying essential analytical properties for numerical methods to capture, every stiffness index does not capture all of the numerically difficult problems. Shampine famously stated ``Indeed, the classic Robertson chemical kinetics problem typically used to illustrate stiffness has two zero eigenvalues. Any definition of stiffness that does not apply to this standard test problem is not very useful.'' \cite{shampine2007stiff}. In addition, the role of the time span is often overlooked. Shampine notes that a system may not seem stiff if one is only solving on the timescale of the fast process: stiffness requires that one is investigating the effects on the time scale of the slow process.

Due to the ambiguity of the definition of stiff systems, we will choose two differential equations studied extensively by the stiff ODE solver literature as representative highly stiff systems \cite{wanner1996solving}: ROBER \cite{robertson1966solution} and POLLU \cite{verwer1994gauss}.  Specifically, the canonical ROBER problem consists of three species and five reactions, and the POLLU problem consists of 20 species and 25 reactions describing the air pollution formation in atmospheric chemistry. The formula of the ROBER and POLLU problem is presented in Section \ref{sec:experiments} and Supplementary Material. We will investigate the difficulty of fitting data sampled from these systems, thus illuminating the unique features of the training process which arise in the discovery of stiff systems and the numerical methods for handling them.

\section{Method}
\label{sec:main_methods}

\subsection{Stabiliy of Adjoints on Stiff ODEs}
\label{subsec:eest_reg}

The core problem of training neural ODEs is calculating the gradient of the ODE's solution. In the method proposed by \citet{chen2018neural}, the ODE is reversed and augmented with the adjoint equations as:
\begin{equation}\label{eq:adjoint}
\frac{d}{dt}\left[\begin{array}{c}
z\\
\omega\\
\frac{dL}{d\theta}
\end{array}\right]=-\left[\begin{array}{ccc}
1 & \omega^{T} & \omega^{T}\end{array}\right]\left[\begin{array}{ccc}
f & \frac{\partial f}{\partial z} & \frac{\partial f}{\partial\theta}\\
0 & 0 & 0\\
0 & 0 & 0
\end{array}\right]
\end{equation}
However, previous work has called into question the stability of reversing the ODE equation, i.e.
\begin{equation}
z^\prime = -f(z,\theta,-t)
\end{equation}
when solving from $t_{1}$ to $t_0$. While theoretically this reverses the ODE exactly, in practice numerical issues can cause errors which propagate into the gradient. \citet{gholami2019anode} noted that if the Lipschitz constant is too large in absolute value, then either the forward pass or reverse solve is numerically unstable. This can be seen by using a standard ODE solver on the linear ODE $\frac{dy}{dt} = -\lambda y$ on $t \in (0,1)$ with $y(0)=1$, $\lambda=100$ imparts about 1\% error while $\lambda=10000$ cannot be reversed numerically in double precision.

\begin{figure}
    \centering
    \includegraphics[width=\linewidth]{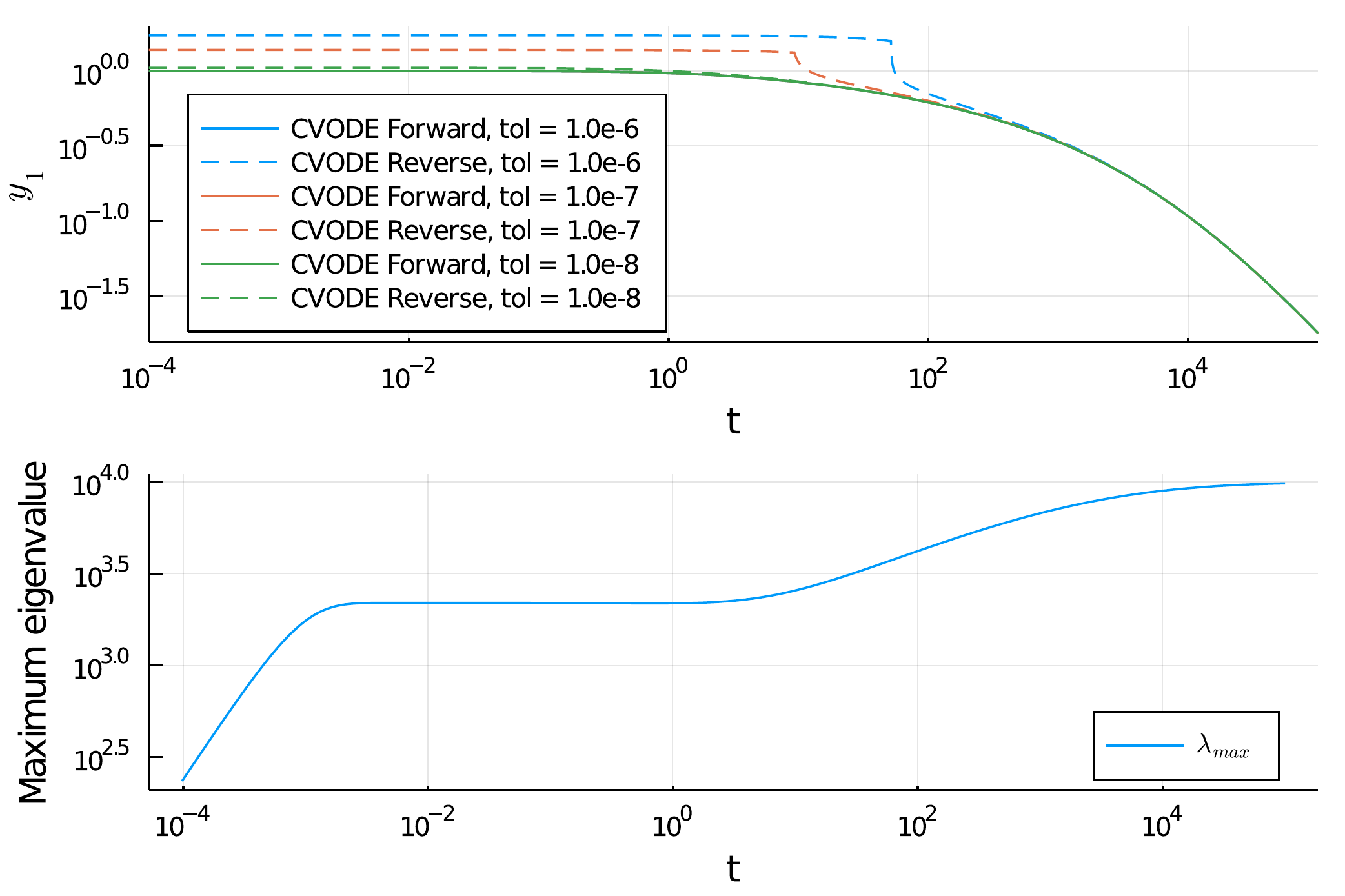}
    \caption{\textbf{Reversability of the ROBER \cite{robertson1966solution} Equation}. \textbf{Top:} Solution forward and backwards of $y_1$ in the ROBER test problem with the CVODE integrator at the shown tolerance using the BDF tableau. \textbf{Bottom:} Maximum eigenvalue of the Jacobian over the solution trajectory.}
    \label{fig:reversability}
\end{figure}

To show how this result extends to stiff ODEs, recall that the Lipschitz constant of a nonlinear function is the infimum of values $\Lambda$ such that $\left\|\frac{f(x)-f(y)}{x-y}\right\| \leq \Lambda$ for some domain $x,y \in \Omega$, which directly implies that the global Lipschitz constant satisfies $\Vert\frac{df}{dy}\Vert_\infty \leq \Lambda$ over the whole domain. Thus because the maximum eigenvalues of the Jacobian give a lower bound on the Lipschitz constant, stiff ODEs fall into the class of problems which are numerically difficult to reverse. Figure~\ref{fig:reversability} demonstrates this on the ROBER stiff ODE test problem \cite{robertson1966solution, ji2020stiff} and demonstrates that even with imperceptible errors in the forward pass we can receive noticeable errors in the reverse solve even when plotted in log scale. Using the well-known \verb CVODE  BDF integrator from Sundials \cite{hindmarsh2005sundials}, we see that with $tol=abstol=reltol=10^{-6}$ we receive 72\% error, $tol=10^{-7}$ gives 38\% error, and $tol=10^{-8}$ gives 5\% error, with 0.1\% error only reached at $tol<10^{-10}$, i.e. below the threshold of single precision arithmetic.

\begin{figure}
    \centering
    \includegraphics[width=\linewidth]{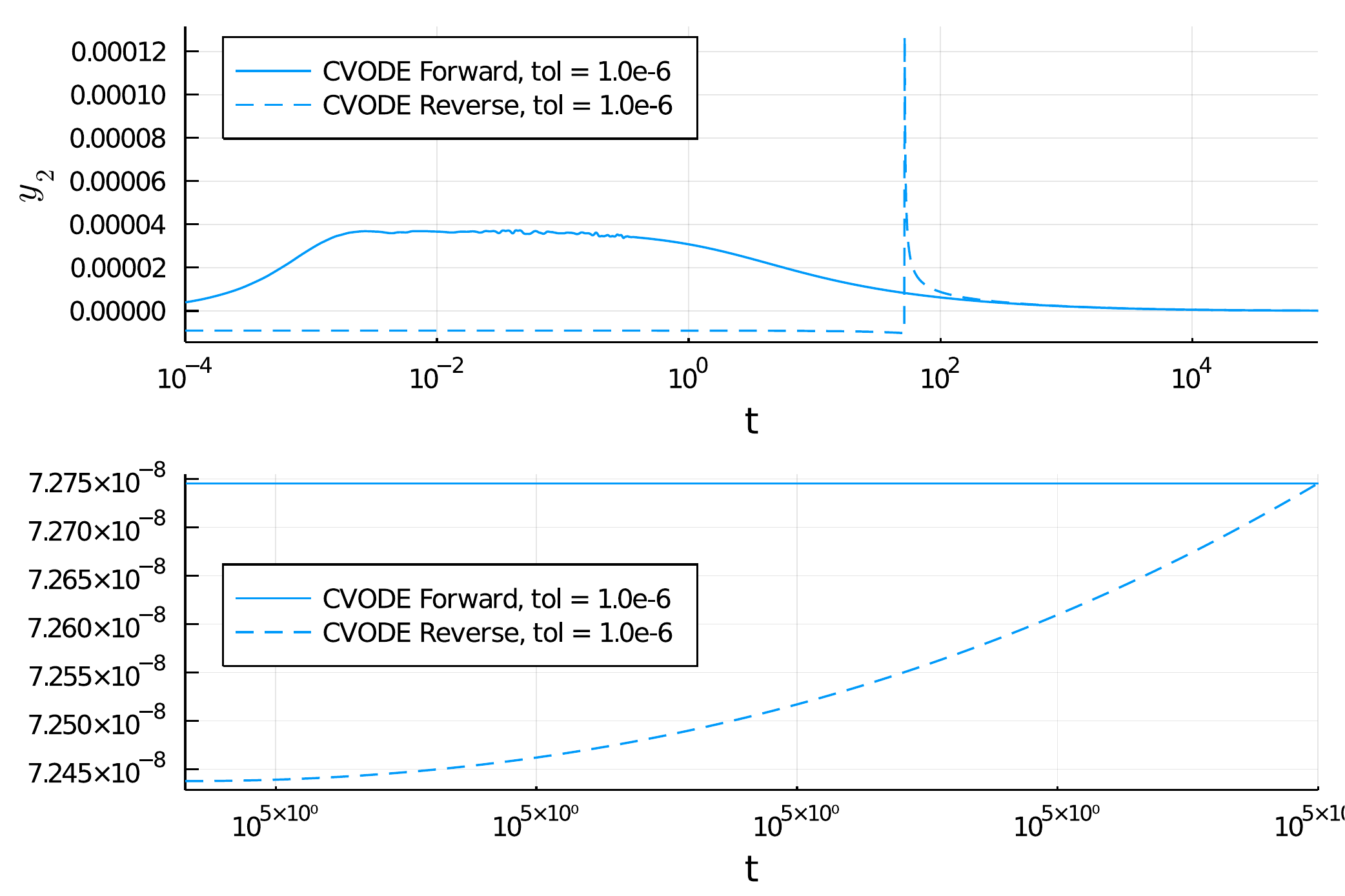}
    \caption{\textbf{Reversability instability inspection on the ROBER Equation}. {\bf Top:} Solution forward and backwards of $y_2$ in the ROBER test problem with the CVODE integrator at the shown tolerance using the BDF tableau. {\bf Bottom:} Zoomed in solution on $t\in (99999.98999,10^5)$.}
    \label{fig:reversability2}
\end{figure}

The adjoint calculation only amplifies this error. With a test loss function $L(\theta) = \sum y_1(t_i)$ at evenly spaced points in log-space, using the adjoint as described caused \verb CVODE  to exit with a Newton divergence near the start of the reverse pass for all $tol=10^{-i}$ for integers $i={6,7,8,9,10,11,12,13,14}$. Manually inspecting the solution process reveals the problem. While the reverse solution only has an approximately $10^{-10}$ error in $y_2$ at $t=99999.98999$ as shown in Figure \ref{fig:reversability2}, the large $k_2$ amplifies this term by $6 \times 10^7$ imparting approximately $\mathcal{O}(10^{-3})$ error into $\frac{d}{dt} \frac{dL}{d\theta}$. This is because while in the original ODE the middle term is of the form $k_2 y_2^2$, in the Jacobian the term is $2 k_2 y_2$ and without the squaring the small error is not diminished. This sends $\frac{dL}{d\theta}$ negative in a way that linearly explodes.

\begin{figure}
    \centering
    \includegraphics[width=\linewidth]{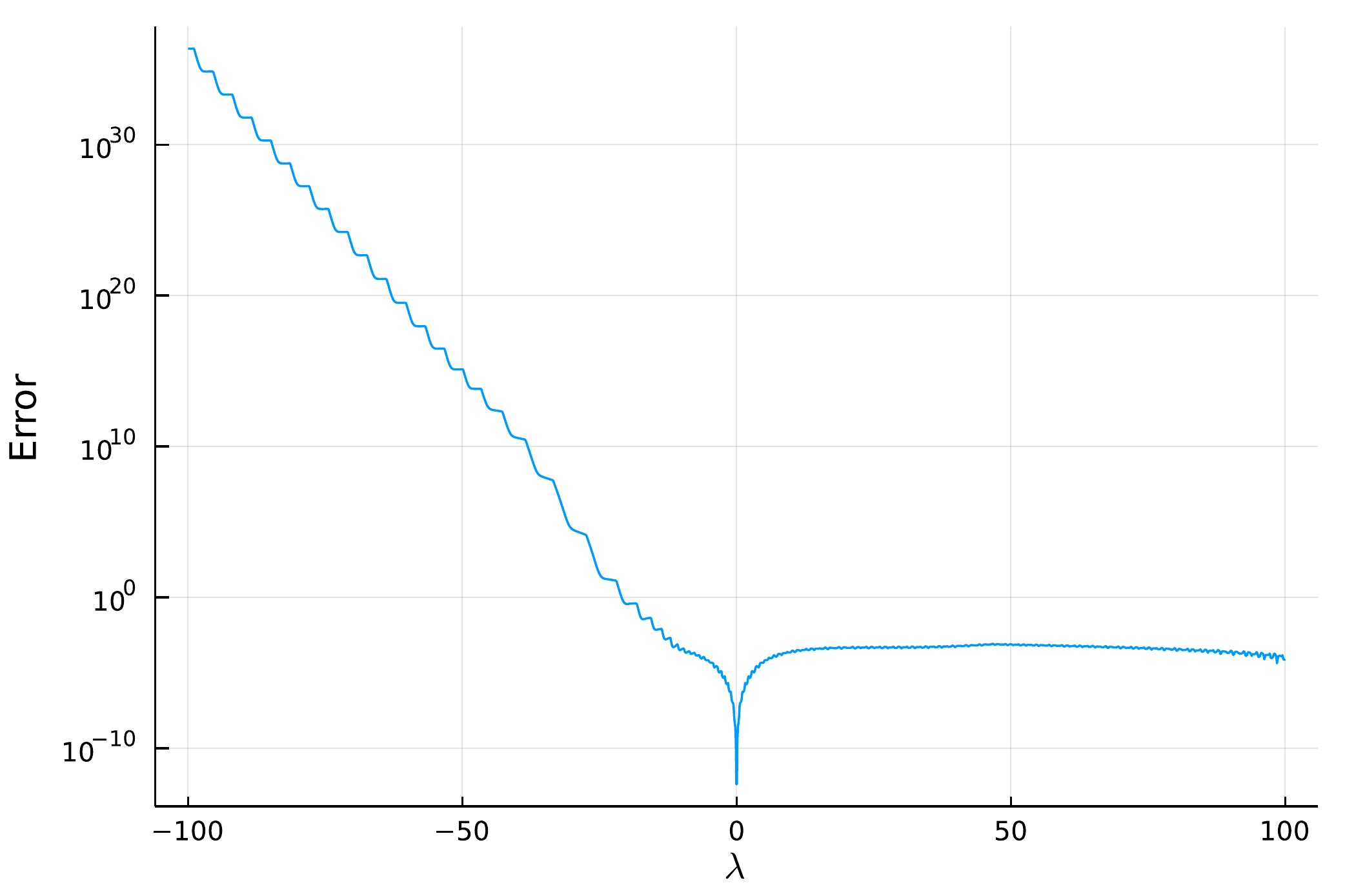}
    \caption{\textbf{Reversability of the linear ODE}. The total error of the solution produced by the Tsit5 solver at time $t=0$ is plotted with varying $\lambda$. Note that the error is plotted in logarithmic scale to illustrate the exponential blow-up behavior.}
    \label{fig:exp_blowup}
\end{figure}

The simple linear ODE $u'(t) = \lambda u(t), u(0)=1$ on the domain $t\in [0, 1]$ is sufficient to demonstrate the blow-up behavior of the reverse solve after the forward solve. In the presence of truncation error and round-off error, the forward solution is $u_f(t) = e^{t\lambda} + \varepsilon_f(t)$, where $\epsilon_f$ denotes the error of the forward solve. It follows that the reverse problem has the initial condition $u_b(1)=e^{\lambda} + \varepsilon_f(1)$. Thus, at $t=0$ the backward solution is $u_b(0)=1+e^{-\lambda}\varepsilon_f(1) + \varepsilon_b(0)$, where $\varepsilon_b$ denotes the error of the backward solve. Therefore, the total error at $t=0$ from solving the linear equation forward and backward is $e^{-\lambda}\varepsilon_f(1) + \varepsilon_b(0)$, and it exhibits inherit exponential blow-up behavior when $\lambda<0$. A stable forward problem implies the exponential divergence of the reverse problem. To demonstrate the exponential diverge, the above simple linear ODE is solved forward and backward, and the error at $t=0$ is plotted with various parameters $\lambda$ in Figure \ref{fig:exp_blowup}.

Avoiding this issue requires one does not attempt to reverse the ODE. To this effect we tested the interpolated checkpointing method of \verb CVODES  \cite{hindmarsh2005sundials} 
(known as \verb InterpolatingAdjoint  in DiffEqFlux)
and discrete adjoint sensitivities via automatic differentiation of the solver 
\cite{rackauckas2020generalized}
which both successfully solved the equation and agreed on all derivatives to at least three decimal places when solved with a $tol=10^{-6}$. We note that one can prove the consistency of the adjoint discretiation in the case of discrete adjoint sensitivities ensuring the stability \cite{zhang2014fatode}.

\subsection{Linear Scaling Adjoints}
\label{subsec:linscale_adj}

Being limited to non-reversing adjoints, we investigate the computational cost of derivative calculations on stiff ODEs. To illustrate the computational cost and motivate new formulations, consider the implicit Euler method:
\begin{equation}
    y_{n+1} = y_n + hf(y_{n+1},t+h)
\end{equation}
which is solved via the root finding problem $G(y_{n+1}) = y_{n+1} - y_n + hf(y_{n+1},t+h) = 0$ for the next step. The instability of fixed point solvers \cite{wanner1996solving} makes it necessary to solve the root finding problem by Newton's method, i.e. $y_{i+1} = y_i - \frac{dG(y_i)}{dy_i}^{-1} G(y_i)$ so that $y_{i} \rightarrow y_{n+1}$. The core computational cost is thus the solving of the linear system $\frac{dG(y_i)}{dy_i} x = G(y_i)$. The condition number of the matrix  $\frac{dG(y_i)}{dy_i}$ is precisely the ratios of the maximum and minimum eigenvalues of the system's Jacobian, meaning the linear system is ill-conditioned for stiff equations. Krylov subspace methods like GMRES have slow convergence when the condition number of the matrix is high \cite{liesen2004convergence} and thus using these methods can require problem-specific preconditioners to make convergence computationally acceptable. For this reason, methods for solving stiff ODEs generally attempt direct methods, either dense or sparse, on the Jacobian. The general cost of performing such a linear solve via LU-factorization is $\mathcal{O}(k^3)$ where $k$ is the number of columns in the Jacobian which is equivalent to the number of states in the ODE.

The cubic cost growth of solving a stiff system changes the efficiency calculus of methods for computing the derivative of the ODE solution and thus for calculating the gradient of the neural ODE's loss function. In the traditional argument, forward-mode automatic differentiation of an ODE is equivalent to solving the forward sensitivity equations:
\begin{align}
    z^\prime &= f(z,p,t)\\
    \frac{d}{dt} \frac{dz}{dp} &= \frac{\partial f}{\partial z} \frac{dz}{dp} + \frac{\partial f}{\partial p}
\end{align}
where there exists one sensitivity $\frac{dz}{dp}$ for each parameter. Thus given $k$ state equations and $m$ parameters, the total size of the ODE system is $\mathcal{O}(km)$. In contrast, the adjoint method works by only solving the original ODE forwards, and solves a system of size $\mathcal{O}(m+k)$ in reverse, greatly reducing the computational cost when the size of the ODE and the parameters are sufficiently large.

However, the total computational cost for implicit methods is based on the cube of the system size, meaning standard forward-mode sensitivities scale as $\mathcal{O}(k^3m^3)$ while the adjoint method of \citet{chen2018neural} or \citet{serban2003cvodes} scales as $\mathcal{O}((k + m)^3)$. Given the large number of parameters in a neural ODE, the cubic scaling adds a tremendous computational cost. To mitigate this issue, we developed 3 separate strategies. For the first, note that many solvers contain a scheme known as dense output for generating a high order interpolation by reusing the internal computations of the steps \cite{shampine2015dense}. This interpolation scheme in many cases requires zero additional $f$ calculations by the ODE solver to construct and thus is computationally efficient. However, because it uses the intermediate steps for the calculation, it requires an increase in the memory by $s$ for a $s$-stage Runge-Kutta or Rosenbrock method\footnote{Note many methods do not need to store all stages. Additionally, some higher order Runge-Kutta methods require a few additional $f$ evaluations in order to generate an interpolation matching the order of the solver and thus can default to interpolants of an order less or more. For a full discussion of dense output, see \cite{shampine2015dense}.}. Using the dense output schemes we can instead calculate Equation \ref{eq:adjoint} by sequentially solving $z^\prime = f(z,\theta,t)$ forwards, then solving $\omega^\prime = -\omega^T \frac{\partial f}{\partial z}$ in reverse (using the dense output of $z(t)$ for any Jacobian computation), and finally  solve $\frac{dL}{d\theta} = \int_{t_{0}}^{t_{1}}\omega^{T}(t)\frac{\partial f}{\partial \theta}dt$ independently using the dense output for $\omega(t)$. Using quadrature techniques like Clenshaw-Curtis can converge exponentially fast, using less steps than the ODE solver. Thus, it potentially further decreases the computational cost of this step of the calculation which is $\mathcal{O}(m)$. But most importantly, this splits the computation into two stiff ODE solves matching the size of $z$ plus a explicit quadrature matching the size of $\theta$. Because of the cubic cost scaling of LU-factorizations in a stiff ODE solve, this trades computation for memory and brings the cost down to $\mathcal{O}(k^3 + m)$. We termed this the \verb QuadratureAdjoint  method.

A core fact which is exploited in that formulation is that the final equation is explicit and thus does not necessarily have similar stability properties to the forward or adjoint equations. It would only be stiff if the Jacobian of $\omega^{T}(t)\frac{\partial f}{\partial \theta}$ is ill-conditioned, which does not necessarily follow from the conditioning of $\frac{\partial f}{\partial z}$. Thus another formulation which could reduce the memory cost while achieving a similar effect would be to use an implicit-explicit (IMEX) ODE solver. IMEX ODE solvers allow for solving the differential equation as $f = f_i + f_e$ where $f_i$ is treated implicitly while $f_e$ is treated explicitly. In this sense, Equation \ref{eq:adjoint} could be split so that:
\begin{equation}\label{eq:imex}
f_i=-\left[\begin{array}{ccc}
\omega^{T} & 0\end{array}\right]\left[\begin{array}{ccc}
\frac{\partial f}{\partial z} & 0\\
0 & 0
\end{array}\right]
\end{equation}

\begin{equation}
f_e=-\left[\begin{array}{ccc}
0 & \omega^{T}\end{array}\right]\left[\begin{array}{ccc}
0 & \frac{\partial f}{\partial\theta}\\
0 & 0
\end{array}\right]
\end{equation}
As an IMEX solver only requires factorizing the Jacobian of the $f_i$ term, specializing the factorization on the block structure reduces the LU-factorization cost to $\mathcal{O}(k^3)$, and factoring in the linear cost of the explicit Runge-Kutta handling we arrive at another $\mathcal{O}(k^3 + m)$ scheme but with reduced peak memory requirements.

These two approaches highlight that the main difficulty is due to the inclusion of $m$ parameter derivatives into the stiff solve. Thus, the last way to achieve a similar effect is to simply do any of the aforementioned schemes but only calculate the derivative of $l \ll m$ parameters at a time. In this case, one would need to repeat the adjoint solve $\frac{m}{l}$ times. However, because each solve has a cubic computational cost with respect to the states, this form of splitting brings the complexity down to $\mathcal{O}(m k^3 l^2)$ for forward sensitivities and $\mathcal{O}(\frac{m}{l} (k+l)^3)$ for adjoint. Note that this technique also applies to direct differentiation via forward and reverse automatic differentiation of the solver. Additionally, the separate differentiation passes can all be solved in parallel, further reducing the real-world compute cost. This is the technique that we used to make the experiments of Section \ref{sec:experiments} computationally viable. 

\subsection{Equation Scaling}
Armed with computationally stable and efficient gradient calculations, we noticed that these techniques were still not enough to stabilize the training process. To further tackle the difficulties in training stiff neural ODEs introduced by scale separation, we propose to normalize the neural network outputs and the loss components for different species. Specifically, we learn the following normalized form of neural ODEs: 

\begin{equation}\label{eq:eq_4_2}
\begin{split}
\frac{dy\left(t \right)}{dt} &= NN \left(y \left(t \right), t \right) \frac{y_{scale}}{t_{scale}}\\
y_{scale} &= y_{max} - y_{min}\\
t_{scale} &= t_{1} - t_{0}
\end{split}
\end{equation}

Where $y_{scale}$ is a vector which consists of the characteristic scales for each species and $t_{scale}$ refers to the characteristic time scale. We further re-balance the loss components by normalizing the loss functions as

\begin{equation}
L(\theta) = MAE \left(\frac{y \left(t \right)^{model}}{y_{scale}}, \frac{y \left(t \right)^{obs}}{y_{scale}} \right) \label{eq:eq4_3}
\end{equation}

\section{Experiments}\label{sec:experiments}

\subsection{ROBER Problem}

Robertson’s equations, denoted as ROBER, are one of the prominent stiff ODEs which model a reaction network with three chemicals \cite{robertson1966solution}. The species concentrations $[y_1,y_2,y_3]$ are governed by Equation \ref{eq:eq_4_1} with reaction rate constants of $k_1=0.04, k_2=3\cdot10^7, k_3=10^4$.

\begin{equation}\label{eq:eq_4_1}
\begin{split}
\frac{dy_{1}}{dt} &= -k_{1}y_{1}+k_{3} y_{2} y_{3}\\
\frac{dy_{2}}{dt} &= k_{1}y_{1} - k_{2} y_{2}^{2} -k_{3} y_{2} y_{3}\\
\frac{dy_{3}}{dt} &= k_{2} y_{2}^{2}
\end{split}
\end{equation}

\subsubsection{Baseline Model}

\begin{figure}
    \centering
    \includegraphics[width=\linewidth]{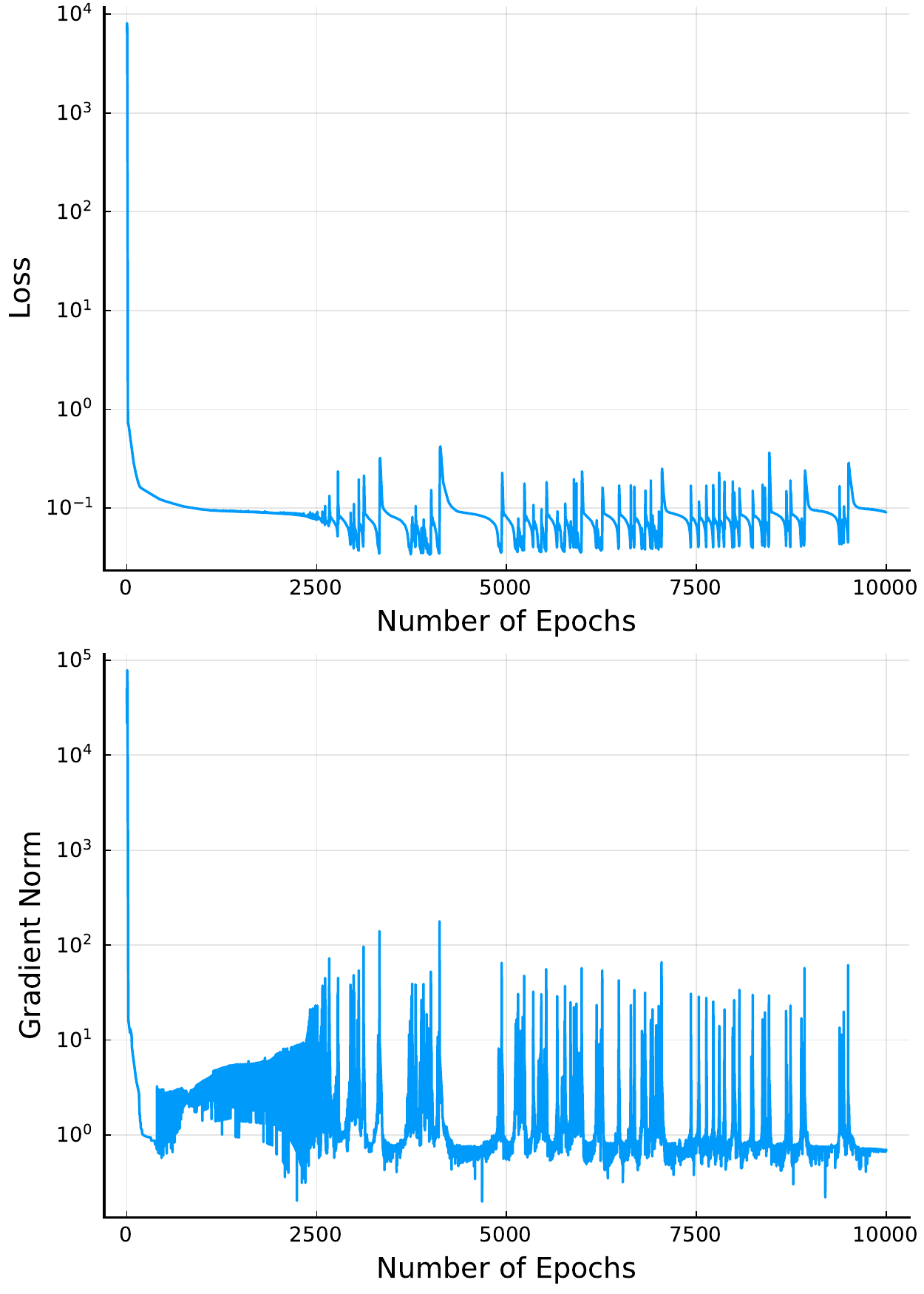}
    \caption{\textbf{Training of the baseline model.} {\bf Top:} Loss profile along epochs {\bf Bottom:} gradient norm  along epochs }
    \label{fig:fig1}
\end{figure}

We first study baseline neural ODEs to elucidate the challenges of training stiff neural ODEs. The baseline model consists of one hidden layer with 50 hidden nodes and a hyperbolic tangent activation function, which is similar to the demo code in DiffEqFlux.jl and torchdiffeq. The stiff ODE integrator \verb Rosenbrock23 , an Order 2/3 L-Stable Rosenbrock-W method, is employed, and we hybrid Rosenbrock23 with Tsit5 via auto-switching to accelerate the integration and training. ForwardDiff.jl is adopted for the backpropagation, which we found is more efficient than adjoint methods in the current case due to the small number of states. The training data is generated by integrating Equation \ref{eq:eq_4_1} with the initial conditions of  $[y_1,y_2,y_3]=[1,0,0]$, time span of $\left[10^{-5},10^{5}\right]$. The 50 sample data points are uniformly spaced in a logarithmic time scale. The neural ODE is trained with the ADAM \cite{kingma2014adam} optimizer with a learning rate of 0.005. The training ended at 10,000 epochs where the loss function flatlines, as shown in Figure \ref{fig:fig1}.

\begin{figure}
    \centering
    \includegraphics[width=\linewidth]{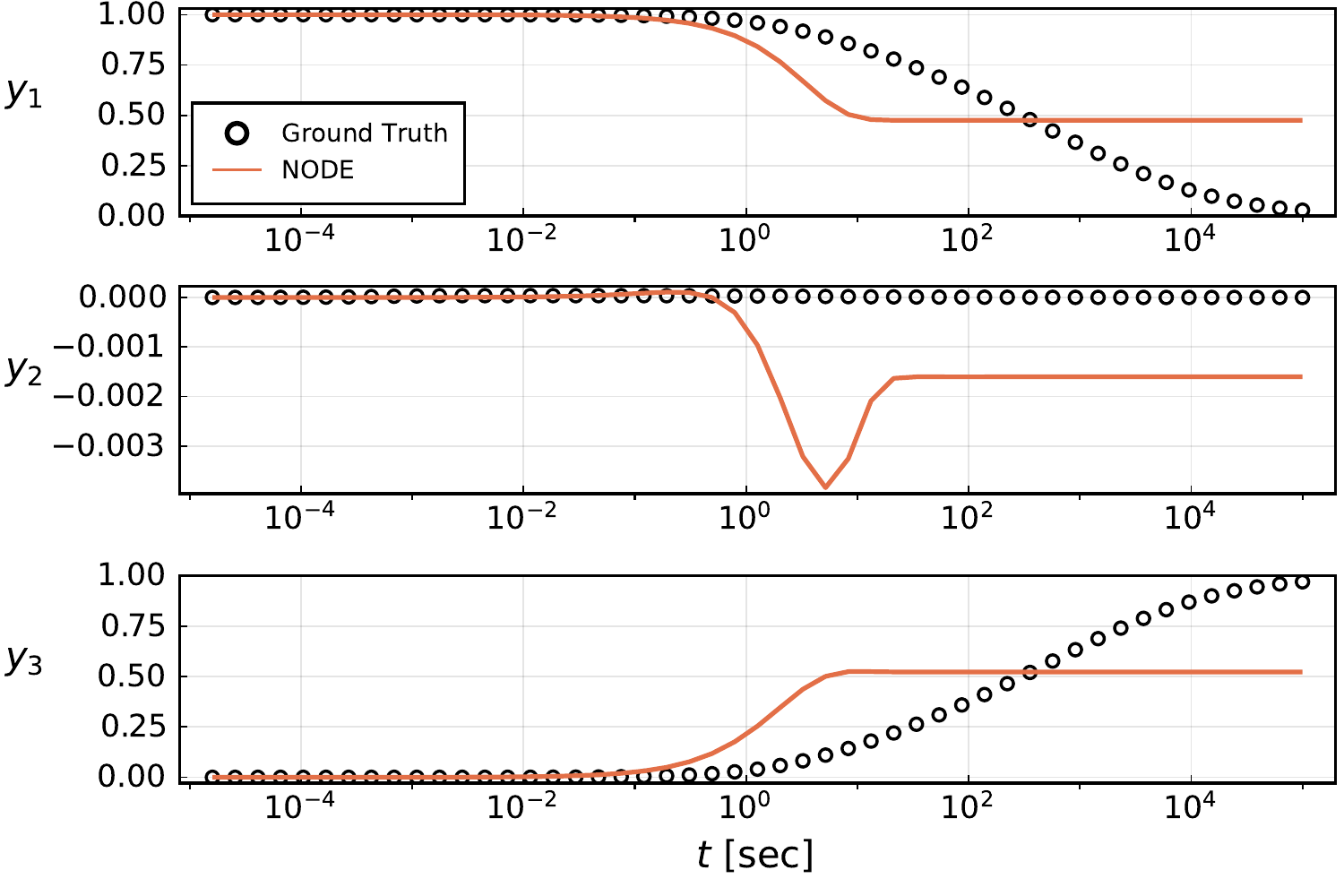}
    \caption{\textbf{The ground truth data (symbols) generated using the ROBER model (solid line) and the predictions of the learned baseline model}}
    \label{fig:fig2}
\end{figure}

Figure \ref{fig:fig2} shows the predictions of the learned baseline model. The baseline model failed to learn the actual dynamics and consequently failed to reproduce the species profiles. In addition, the gradient norm highly fluctuates and training losses do not decrease which implies training instability (Figure \ref{fig:fig1}). We hypothesize that the failure is attributed to the time scale separations and stiffness induced scale separations in the neural network outputs as well as the imbalance of different loss function components. Specifically, after the initial short induction period, the changes of the $y_1$ and $y_3$ are in the order of unity while $y_2$ in the order of $10^{-5}$, which implies that one has to approximate such widely separated scales in a single neural network. In addition, since the magnitude of $y_2$ is 5 orders of magnitude smaller than $y_1$ and $y_3$, such imbalance in the three loss components could lead to gradient pathologies, which was anticipated that stochastic gradient descent is unstable in the presence of such gradient pathologies \cite{wang2020understanding}.

\subsubsection{Mitigation of Gradient Pathologies via Scaling}

\begin{figure}
    \centering
    \includegraphics[width=\linewidth]{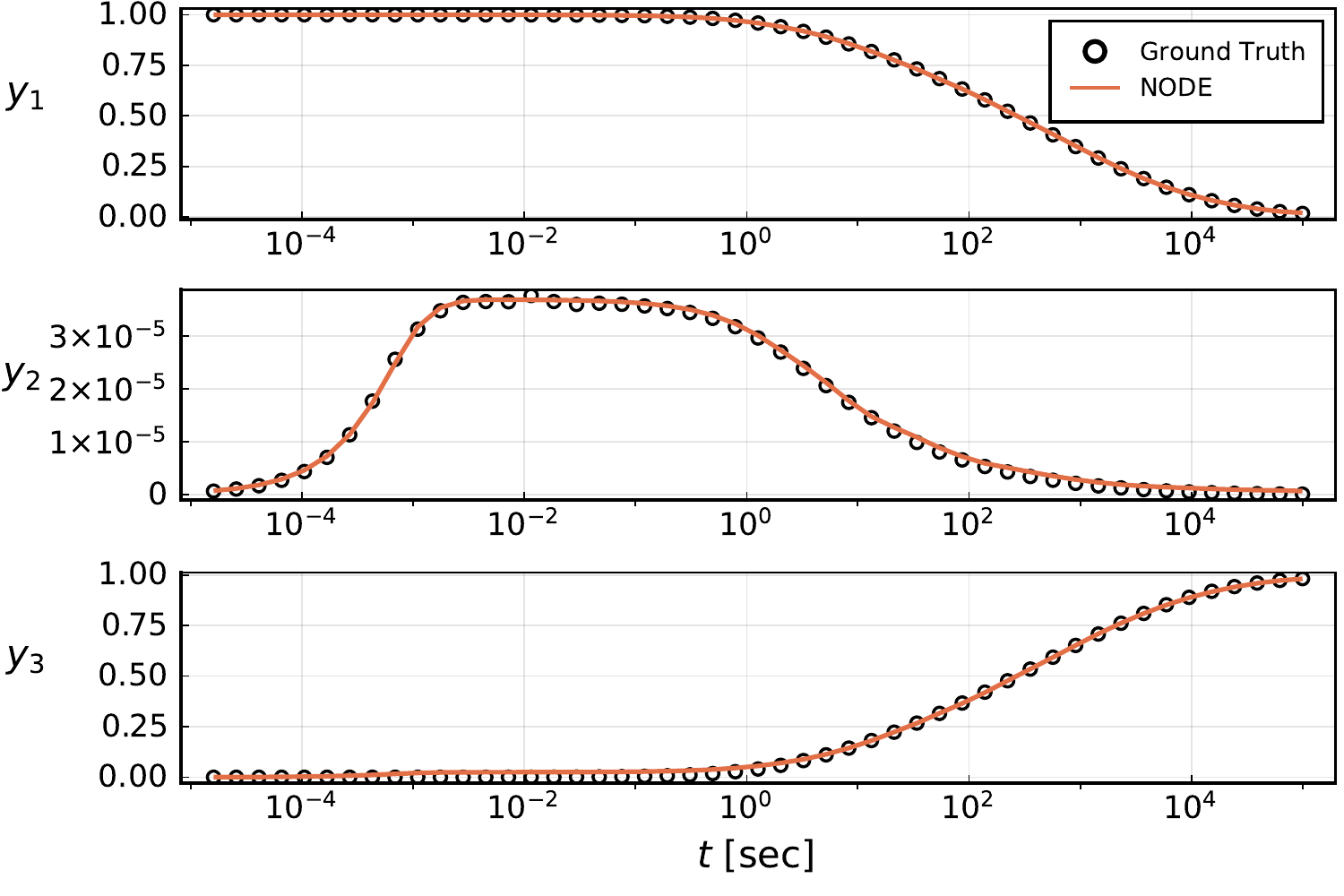}
    \caption{\textbf{Ground truth of Robertson equations (circle mark) and results of the scaled neural ordinary differential equation with six hidden layers of 5 nodes and GELU activation function (orange line)}}
    \label{fig:fig3}
\end{figure}

\begin{figure}
    \centering
    \includegraphics[width=\linewidth]{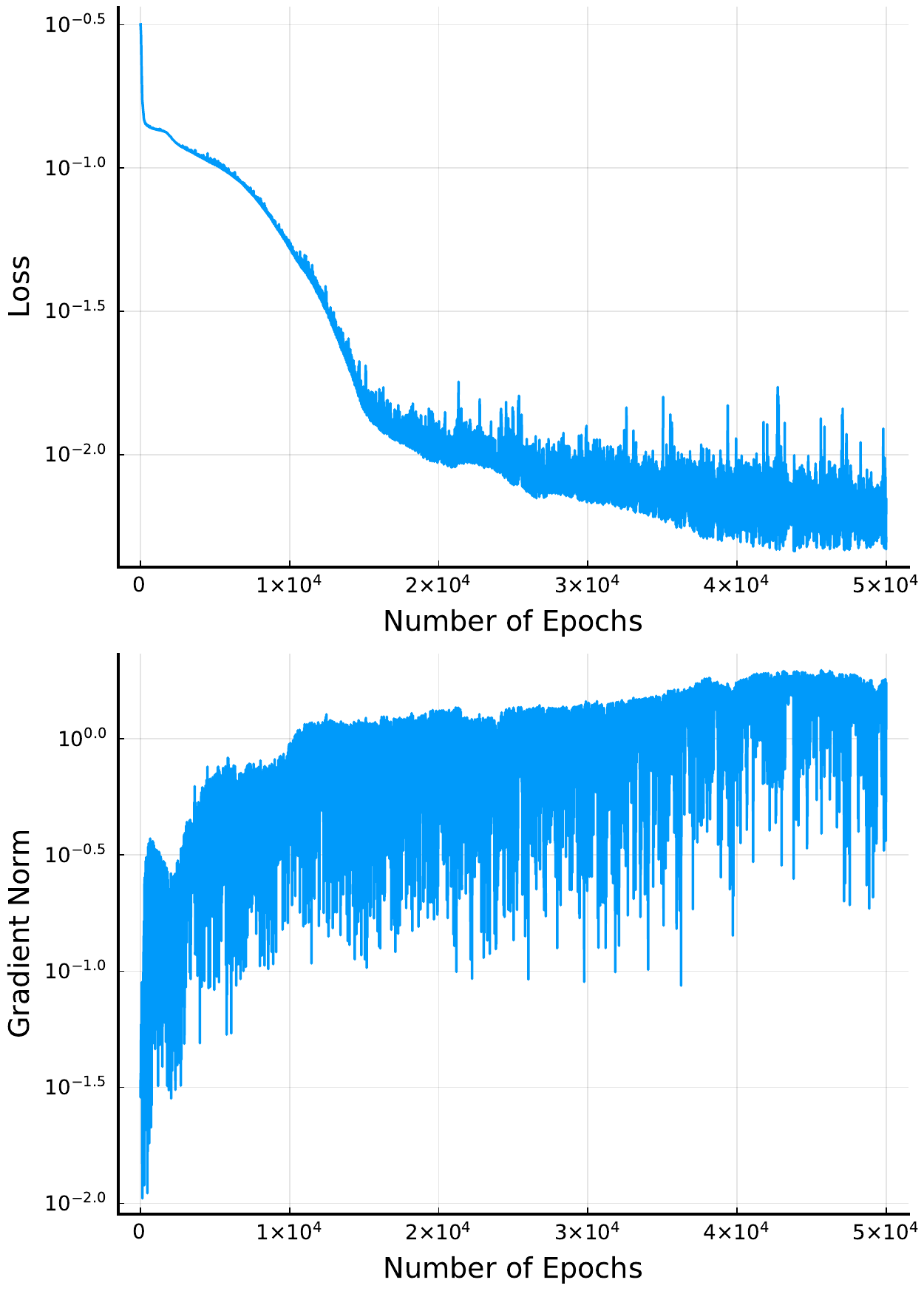}
    \caption{\textbf{Training of a scaled neural ODE.} {\bf Top:} Loss profile along epochs {\bf Bottom:} gradient norm  along epochs }
    \label{fig:fig4}
\end{figure}

Figure \ref{fig:fig3} shows the results of a successful learning of stiff ROBER problems with equation scaling. The neural networks consist of 6 hidden layers and 5 nodes per layer, with an activation function of GELU \cite{hendrycks2016gaussian}. GELU was chosen to mitigate potential issues with vanishing gradients typically seen in recurrant neural networks with saturating activation functions \cite{pascanu2013difficulty}. Figure \ref{fig:fig4} shows the history of loss functions and the $L_2$ norm of the gradient. The loss functions steadily decrease, and the gradient norm is stable during the training. The training in a workstation using a single CPU processor took 2 hours.

We carried out sensitivity analyses to better understand the importance of proposed techniques, i.e., the equation scaling and the deep networks, on the training of stiff ODEs. Our test shows that scaling is essential for training stiff neural ODEs, as the training with the same neural network structure but without equation scaling failed to learn the ROBER problem. Regarding the neural network structure, we found that shallow neural networks with the hyperbolic tangent activation function (tanh) can also successfully learn the ROBER problem but takes a longer training time, which could be attributed to the potential gradient issues. The results of the training without equation scaling and the training with the shallow network are provided in the supplemental material.

We can also accelerate the training with temporal regularization \cite{ghosh2020steer} and annealing of learning rates although such techniques were not required for successful learning of stiff neural ODEs. The temporal regularization, STEER \cite{ghosh2020steer}, proceeds as randomly truncate the training to an integration time ahead of the entire time-space, and in such a way introduce stochastic into the optimization to accelerate the training. The annealing of a learning rate is a widely adopted technique for accelerating the training of neural networks at later stages. 

\subsection{POLLU Problem}

To show the generality of the technique, we demonstrate the effectiveness of the proposed methods in the POLLU problem, which is more complex than the ROBER problem and consists of 20 species and 25 reactions. The POLLU is an air pollution model developed at The Dutch National Institute of Public Health and Environmental Protection. It can be described mathematically by the following 20 non-linear ODEs shown in Equation \ref{eq:eq4_4}, and full details about the model can be found in \cite{verwer1994gauss} and supplemental materials.
 
\begin{equation}
\begin{split}
\frac{dy\left(t \right)}{dt} &= f \left(y \left(t \right) \right)\\ y\left(0\right)=y_{0},  y&\in\mathbb{R}^{20},  0\leq t\leq 60 \label{eq:eq4_4}
\end{split}
\end{equation}

\begin{figure*}
    \centering
    \includegraphics[width=\linewidth]{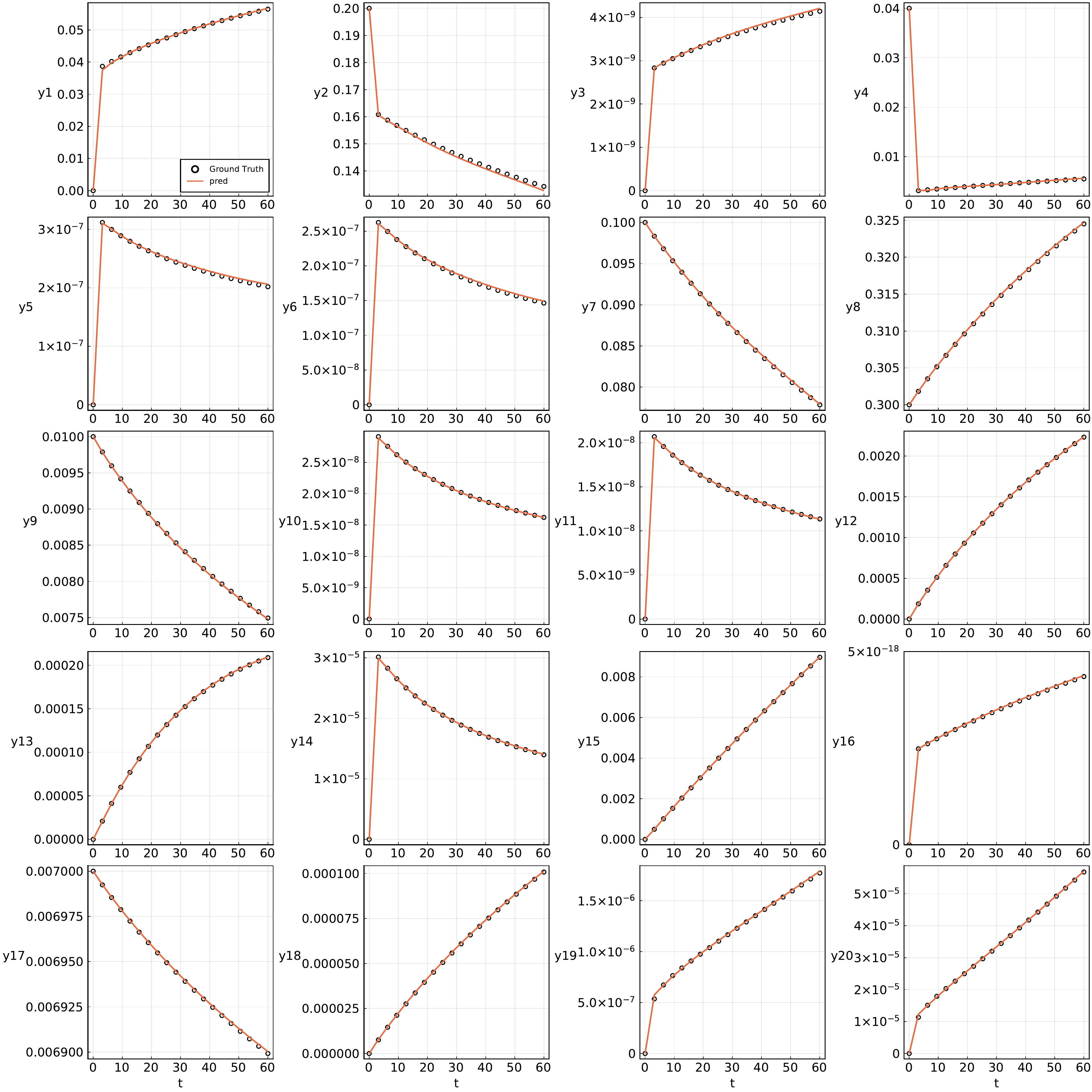}
    \caption{\textbf{Ground truth of an air pollution model POLLU (symbols) and results of the scaled neural ordinary differential equation with 3 hidden layers of 10 nodes and GELU activation function (solid line)}}
    \label{fig:fig5}
\end{figure*}

The neural ODEs consist of 3 hidden layers and 10 nodes per hidden layer. The label data is sampled uniformly in $[0, 60]$ on a linear scale. The rest of the hyper-parameter settings are the same as the ones for the ROBER problem. Figure \ref{fig:fig5} presents the label data generated with Equation \ref{eq:eq4_4} and the predictions of the learned neural ODEs. Note that the baseline model again failed for the POLLU problem similarly to the ROBER problem, and it is presented in the supplemental material. The trained neural ODEs can capture both the dynamics during the induction period and the near-equilibrium states. It also well captures both the dynamics of the major species, such as $y_{1},y_{2}$, with relatively large concentration and slow time scales, and the dynamics of transient intermediate species, such as $y_{3},y_{10}$, with relatively low concentrations and fast time scales. Therefore, the demonstration of the POLLU problem further shows the robustness of proposed approaches in learning stiff neural ODEs where previous techniques had failed.

\section{Conclusion and Discussion}

This work proposed new derivative calculation techniques for the stability of stiff ODE systems while reducing the computational complexity. It is coupled with equation scaling and deeper neural networks with rectified activation functions to learn neural ODEs in stiff systems. The effectiveness of the proposed approaches is demonstrated in two classical stiff systems: the ROBER and the POLLU. While the baseline models even failed to learn the overall trend of the dynamics, the proposed approaches accurately reproduced all of the species profiles including both species with slow and fast time scales.

While we have shown the effectiveness of the proposed techniques, there is a demand for further theoretical analysis of the challenges of learning stiff neural ODE, such as how stiffness affects the gradient dynamics during training and leads to training failures. In addition, more theoretical analysis on equation scaling relations to the gradient flow of neural ODEs would better elucidate why it was important in the training process. It is also interesting to explore other equation scaling and neural network scaling techniques, such as scaling the neural network inputs and batch-normalization \cite{ioffe2015batch}. Given the first order optimizers are solving an ordinary differential equation on the gradient flow themselves, one likely could define the stiffness index on the Jacobian of the gradient of the optimization, which is the Hessian with respect to the parameters. Handling the optimization itself with stiffly-aware methods could be beneficial in cases where gradient pathologies are not mitigated.

Importantly, this manuscript demonstrates the impact of errors in the continuous adjoint handling on the training of a neural ODE. Such an analysis could be required in other ``implicit layer'' methods as well. For example, we referenced how a nonlinear solver destabilizes when its Jacobian exhibits ill-conditioning, equivalent to stiffness in the ODE system when viewing it as a solver to steady state \cite{wanner1996solving}. This would suggest that on highly stiff systems, methods like the DEQ \cite{bai2019deep} similar behavior would apply to the forward pass. But importantly, similar instabilities in the adjoint problem likely occur due to the reliance on the adjoint's assumption that $G(x)=0$ exactly, which is more difficult to satisfy as the condition number rises. Likely similar stabilization may be required for this case. Similarly, differentiable optimization as a neural network layer \cite{amos2017optnet,amos2019differentiable} requires satisfaction of the KKT equations for its adjoint, which are likely to be less satisfied when the Hessian of the optimization problem is ill-conditioned, therefore likely requires a similar stability analysis.

We end by noting that such developments in handling highly stiff equations can allow for machine learning architectures that satisfy arbitrary constraint equations during their evolution. \citet{rackauckas2020universal} proposed using differential-algebraic equations (DAEs) in mass-matrix form for this purpose, i.e.
\begin{equation}
    Mz^\prime = f(z,p,t)
\end{equation}
where $M$ is singular. For example, it's well-known that the Rosenbrock equation can be written in its DAE form:
\begin{align}
    \frac{dy_{1}}{dt} &= -k_{1}y_{1}+k_{3} y_{2} y_{3}\\
    \frac{dy_{2}}{dt} &= k_{1}y_{1} - k_{2} y_{2}^{2} -k_{3} y_{2} y_{3}\\
    0 &= y_1 + y_2 + y_3 - 1
\end{align}
where the mass matrix is of the $M = [1 0 0;0 1 0; 0 0 0]$. Similarly, an arbitrary system constrained to have dynamical states sum to 1 could be imposed by:
\begin{align}
    \frac{d[y_1,y_2]}{dt} &= NN(y)\\
    0 &= y_1 + y_2 + y_3 - 1
\end{align}
and more generally
\begin{align}
    \frac{dy}{dt} &= NN([y,z])\\
    0 &= g(y,z,t)
\end{align}
constraints equation can be set or fit using domain information as necessary. While currently able to be specified and trained in differentiable DAE solver systems,
we noticed fitting instabilities similar to stiff neural ODEs arise in such cases. Indeed, \citet{wanner1996solving} notes that DAEs are a form of infinite stiffness, being posed as the limit of singularly perturbed stiff ODEs. \citet{petzold1982differential} famously showcased how ``DAEs are not ODEs'', showing how they exhibit additional behavior that must be appropriately treated, and thus an analysis of such systems will be required for fully stabilizing neural network DAE formulations including the relationship of training techniques to differential index.

\section*{Supplementary Material}
See \href{run:./Stiff_Neural_ODES_Supp.pdf}{supplementary material} for the sensitivity analysis to the neural network structures and the details in POLLU equations \cite{verwer1994gauss}.

\begin{acknowledgements}
The information, data, or work presented herein was funded in part by the Advanced Research Projects Agency-Energy (ARPA-E), U.S. Department of Energy, under Award Number DE-AR0001222. The views and opinions of authors expressed herein do not necessarily state or reflect those of the United States Government or any agency thereof.
\end{acknowledgements}

\section*{Data Availability}

The data that support the findings of this study are openly available on GitHub at https://github.com/DENG-MIT/StiffNeuralODE.

\section*{References}
\bibliography{aipsamp}

\end{document}